\documentclass[11pt]{amsart}

\newtheorem{thm}{Theorem}[section]

\newtheorem{prop}[thm]{Proposition}

\theoremstyle{plain}
\newtheorem{main}{Theorem}
\newtheorem{corm}[main]{Corollary}

\theoremstyle{definition}
\newtheorem*{question}{Question}

\newtheorem*{rem}{Remark}

\usepackage{amssymb,latexsym,amscd,amsmath,amsthm}
\usepackage[all]{xy}
\newcommand{\R}{\mathbb{R}}
\newcommand{\Z}{{\mathbb{Z}}}

\newcommand{\sph}{{\mathbb{S}}}

\newcommand{\f}{\mathbb{F}}

\renewcommand{\H}{\operatorname{H}}
\newcommand{\syp}{\operatorname{Sp}}
\newcommand{\SO}{\operatorname{SO}}

\renewcommand{\lim}[1]{\mathop{\underset{#1} {\underset \longleftarrow
{\text{\rm lim}}}}}

\newcommand{\comp}{{\text{ $\scriptstyle \circ$ }}}

\newcommand{\sq}{\operatorname{Sq}}

\begin{document}

\title{On complexes equivalent to $\sph^3$-bundles over $\sph^4$}

\author{Nitu Kitchloo}
\author{Krishnan Shankar}

\address{Department of Mathematics\\ Northwestern University\\
Evanston\\ IL 60208}
\email{nitu@math.nwu.edu}

\address{Department of Mathematics\\ University of Michigan\\ Ann
	 Arbor\\ MI 48109}
\email{shankar@math.lsa.umich.edu}


\maketitle

\normalsize
\thispagestyle{empty}
\section*{Introduction}

$\sph^3$-bundles over $\sph^4$ have played an important role in
topology and geometry since Milnor showed that the total spaces of
such bundles with Euler class $\pm 1$ are manifolds homeomorphic to
$\sph^7$ but not always diffeomorphic to it. In 1974, Gromoll and
Meyer exhibited one of these spheres (a generator in the group of
homotopy 7-spheres) as a double coset manifold i.e.\ a quotient of
$\syp(2)$ hence showing that it admits a metric of nonnegative
curvature (cf.\ \cite{grommeyer}). Until recently, this was the
only exotic sphere known to admit a metric of nonnegative sectional
curvature. Then in \cite{grovziller}, K. Grove and W. Ziller constructed
metrics of nonnegative curvature on the total space of
$\sph^3$-bundles over $\sph^4$. They also asked for a classification
of these bundles up to homotopy equivalence, homeomorphism and
diffeomorphism. These questions have been addressed in many papers
such as \cite{tamura}, \cite{sasao}, \cite{wilkens} and more
recently in \cite{be}. In this paper we attempt to
fill the gap in the previous papers; we consider the problem of
determining when a given CW complex is homotopy equivalent to such a
bundle. The problem was motivated by \cite{grovziller}: the Berger
space, $\syp(2)/\syp(1)$, is a 7-manifold that has the cohomology ring
of an $\sph^3$-bundle over $\sph^4$, but does it admit the structure
of such a bundle? The fact that it cannot be a principal
$\sph^3$-bundle over $\sph^4$ is straightforward and is proved in
\cite{grovziller}.
\medskip

Let $X$ be a simply connected CW complex with integral cohomology
groups given by
\begin{equation}
	\begin{aligned}
	H^i(X) &= \Z \quad \text{if  }i=0,7\\
	       &= \Z_n \quad \text{if  } i=4
	\end{aligned}
\end{equation}
where $n$ is some fixed integer. We say that $X$ is oriented with
fundamental class $[X]$ if a generator $[X] \in H_7 (X)$ is specified. For 
oriented $X$ we define the
\textit{linking form} as,
\begin{align*}
	b: &H^4(X)\otimes H^4(X) \rightarrow \Z_n \\
	   &x\otimes y \longmapsto \langle
		\beta^{-1}(x),y\cap [X] \rangle
\end{align*}
where $\beta: H^3(X,\Z_n)\rightarrow H^4(X)$ is the Bockstein
isomorphism. The homomorphism $\cap [X]: H^4(X)\rightarrow H_3(X)$ is
obtained by capping with the fundamental class and $\langle,\rangle :
H^3(X,\Z_n)\otimes H_3(X) \rightarrow \Z_n$ is the Kronecker
pairing. We now state our main theorem:

\begin{main} Let $X$ be a simply connected CW complex as above. Then
$X$ is homotopy equivalent to an $\sph^3$-bundle over $\sph^4$ if and
only if the following two conditions hold:
\begin{itemize}

\item[(I)] The secondary cohomology operation $\Theta$ is trivial, where
$$
	\Theta: H^4(X,\f_2) \rightarrow H^7(X,\f_2)
$$
corresponds to the relation $\sq^2 \sq^2 = \sq^3 \sq^1$ in the mod 2 
Steenrod algebra.

\item[(II)] The linking form $b: H^4(X)\otimes H^4(X) \rightarrow
\Z_n$ is equivalent to a standard form for some choice of orientation
on $X$ i.e.\ there exists an isomorphism $\psi: \Z_n \rightarrow
H^4(X)$ such that $b(\psi(x),\psi(y)) = xy$.
\end{itemize}
\end{main}

Using the method outlined on page 32 of \cite{mm}, it is easy to show
that if $X$ smoothable, then $\Theta$ is trivial and we have,

\begin{corm} Let $M$ be a simply connected 7-manifold with integral
cohomology groups given by
\begin{align*}
	H^i(M) &= \Z \quad \text{if  } i=0,7 \\
	       &= \Z_n \quad \text{if  } i=4.
\end{align*}
Then $M$ is homotopy equivalent to an $\sph^3$-bundle over $\sph^4$ if
and only if its linking form is equivalent to a standard form for some
choice of orientation on $M$.\end{corm}

The previous corollary can be strengthened in some cases using the results of D. Wilkens. In his paper \cite{wilkens}, Wilkens classified simply connected manifolds, $M^7$, with integral cohomology as above up to the
addition of homotopy 7-spheres (and hence up to $PL$-homeomorphism
type). From the discussion in the Appendix, we have:

\begin{main} Let $M$ be a simply connected manifold with integral
cohomology as in (1) and $H^4(M)\cong \Z_n$. If either $n$ is odd or
$\frac{n}{2}$ is odd, then $M$ is $PL$-homeomorphic to an
$\sph^3$-bundle over $\sph^4$ if and only if its linking form is
equivalent to a standard form.\end{main}

The above result can be further strengthened to include all integers $n$ using
results in a forthcoming paper by B. Botvinnik and C. Escher (cf.\
\cite{be}). The seven dimensional Berger manifold is a curious space. It is
described as the homogeneous space, $M=\syp(2)/\syp(1) =
\SO(5)/\SO(3)$, where the embedding of $\syp(1)$ in $\syp(2)$ is
maximal. It is an isotropy irreducible space and has the cohomology
ring as in (1) with $H^4(M)=\Z_{10}$ (see Section 4). It admits a
normal homogeneous metric of positive sectional curvature (cf.\
\cite{berger}). In \cite{grovziller} the following question was asked:

\begin{question} Does the Berger space, $M=\syp(2)/\syp(1)$, admit the
structure of an $\sph^3$-bundle over $\sph^4$?\end{question}

In Section 4 it is shown that the linking form for the Berger space is
equivalent to a standard form. So applying Corollary 2 we
see that the Berger space is indeed homotopy equivalent to an
$\sph^3$-bundle over $\sph^4$. Since $|H^4|=10$ for this space, we
apply Theorem 3 to get,

\begin{corm} The Berger space, $M=\syp(2)/\syp(1)$, is
$PL$-homeomorphic to an $\sph^3$-bundle over $\sph^4$.\end{corm}

It remains open whether the Berger space is in fact diffeomorphic to
such a bundle. This involves computing the Eells-Kuiper invariant,
$\mu$, for this manifold (cf.\ \cite{ek}). The $\mu$
invariant for a 7-manifold is computed by exhibiting the space as a
spin boundary; we are unable to do this for the Berger space.
\medskip

Another application of Theorem 2 is the case when $n = |H^4(X)|=p^m$
where $p$ is a prime of the form $p=4k+3$. Since $-1$ is not a square
in the ring $\Z_{p^m}$, any non-degenerate form, $\alpha:
\Z_{p^m}\otimes \Z_{p^m} \rightarrow \Z_{p^m}$, is equivalent to a
standard form (up to sign). This shows:

\begin{main} Let $X$ be a Poincar\'e duality complex with integral
cohomology as in (1) and $H^4(X) = \Z_{p^m}$ where $p$ is a prime of
the form $p=4k+3$. Then $X$ is homotopy equivalent to an
$\sph^3$-bundle over $\sph^4$. \end{main}

The proof of Theorem 1 is organized as follows: In Section 1 we show
that the conditions (I) and (II) are necessary for $X$ to support the
structure of an $\sph^3$ fibration over $\sph^4$. In Section 2 we
establish sufficiency of the conditions and exhibit $X$ as the total
space of an $\sph^3$ fibration over $\sph^4$. In Section 3 we prove
that any $\sph^3$ fibration over $\sph^4$ is equivalent to a linear
$\sph^3$-bundle over $\sph^4$. In Section 4 we calculate the
cohomology of the Berger space and show that its linking form is
equivalent to a standard form. Finally in the Appendix we discuss the
results of Wilkens, Sasao and James--Whitehead that allow us to prove
Theorem 3.
\medskip

It is a pleasure to thank Wolfgang Ziller for many interesting and
insightful discussions.

\section{Necessity of conditions in Theorem 1}

Let $X$ be a simply connected CW complex with integral cohomology as in
(1). It follows that the 4-skeleton of $X$ can be chosen to be
equivalent to the space $P^4(n)$ defined as the cofiber of the self
map of degree $n$ on $\sph^3$. The space $X$ is then equivalent to the
cofiber of some map $f$ given by,
$$
	\begin{CD}
	\sph^6 @>{f}>> P^4(n) @>>> X
	\end{CD}
$$
Assume now that $X$ supports the structure of the total space in a
fibration:
\begin{equation}
	\begin{CD}
	\sph^3 @>>> X @>>> \sph^4
	\end{CD}
\end{equation}

Consider the commutative diagram:
\begin{equation*}
	\xymatrix{ {\sph^6} \ar[r]^{f} \ar[dr]_{0} & {P^4(n)}
	\ar[r]^{i} \ar[d]_{g} & {X} \ar[dl]^{\pi} \\ & {\sph^4} & }
\end{equation*}
where $\pi : X \rightarrow \sph^4$ is the projection map and $g=\pi
\comp i$. An easy argument using the Serre spectral sequence for the
fibration (2) shows that $g^\ast :H^4(\sph^4) \rightarrow H^4(P^4(n))$
is an epimorphism. From the theory of secondary cohomology operations
(cf.\ \cite{mt}) it follows that the secondary operation $\Theta$ is
trivial where,
$$
	\Theta: H^4(X,\f_2) \rightarrow H^7(X,\f_2)
$$
corresponds to the relation $\sq^2 \sq^2 = \sq^3 \sq^1$ in the mod 2
Steenrod algebra.
\medskip

It remains to check condition (II) in the theorem. Consider the
Serre spectral sequence for the fibration (2) converging to
$H^\ast(X)$. We have:
\begin{equation}
	\begin{aligned}
	E^{p,q}_2 &=E^{p,q}_4 = H^p(\sph^4)\otimes H^q(\sph^3)\\
		&d_4(y_3) = n y_4
	\end{aligned}
\end{equation}
where $y_3$ and $y_4$ are suitably chosen generators of $H^3(\sph^3)$
and $H^4(\sph^4)$ respectively. Similarly, for the spectral sequence
in $\Z_n$-coefficients converging to $H^\ast(X,\Z_n)$, we have:
$$
	E^{p,q}_2 =E^{p,q}_\infty = H^p(\sph^4,\Z_n)\otimes H^q(\sph^3,\Z_n)
$$
Note that in both spectral sequences there are no extension problems,
hence we may identify $H^\ast(X)$ or $H^\ast(X,\Z_n)$ with the fifth
stage in the respective spectral sequences.

Let $[X] \in H_7(X)$ be the orientation class defined by
\begin{equation}
	\langle y_4\otimes y_3,[X] \rangle = 1
\end{equation}
It follows from (3) that the Bockstein isomorphism, $\beta:
H^3(X,\Z_n) \rightarrow H^4(X)$ is given by
\begin{equation}
\beta([y_3])=y_4
\end{equation}
where we henceforth use the notation $[y]$ to denote the $\!\mod(n)$
reduction of an integral class $y$. Using (4) and (5), we get
$$
\langle \beta^{-1}(y_4), y_4\cap [X]\rangle = \langle
[y_3],y_4\cap [X]\rangle = \langle [y_3]\cup [y_4],[X]\rangle \equiv
1\!\mod(n)
$$
which is simply the statement that the linking form is equivalent to a
standard form.

\section{Construction of a $\sph^3$ fibration}

The purpose of this section is to show that any CW complex $X$
satisfying conditions (I) and (II) of Theorem 1 is equivalent to
the total space of a fibration with base $\sph^4$ and fiber
equivalent to $\sph^3$.
\medskip

Given such a complex, recall that $X$ fits into a cofiber sequence
$$
	\begin{CD}
	\sph^6 @>{f}>> P^4(n) @>>> X
	\end{CD}
$$
Let $p: P^4(n)\rightarrow \sph^4$ be any map inducing an epimorphism
in cohomology. Such a map always exists since $P^4(n)$ is a four dimensional complex. Condition (I) ensures that the composite, $p\comp f$, 
is null homotopic. Thus we get an extension $\tilde{\pi}$ making the 
following diagram commute.
$$
	\xymatrix{
	{\sph^6} \ar[r]^{f} &{P^4(n)} \ar[d]^{p} \ar[r]
	&{X}\ar[dl]^{\tilde{\pi}} \\
	 &{\sph^4} &
	}
$$
We can further extend the above diagram to:
\begin{equation}
\begin{aligned}
	\xymatrix{
	{X} \ar[d]_{\tilde{\pi}} \ar[r] &{\sph^7}
	\ar[d]^{f(\tilde{\pi})} \\
	{\sph^4} \ar[r]^{[n]} & {\sph^4}
	}
\end{aligned}
\end{equation}
where the upper horizontal map has degree 1 and $[n]$ denotes the
self map of degree $n$. Now the extension $\tilde{\pi}$ is not unique;
the set of extensions admits a transitive action of the group
$\pi_7(\sph^4)$, which we describe below. Given $\alpha \in
\pi_7(\sph^4)$ we define $\alpha \tilde{\pi}$ by
$$
	\xymatrix{
	{X} \ar[d]_{\alpha \tilde{\pi}} \ar[r]^{\text{pinch}} & {\sph^7
	\vee X} \ar[dl]^{\alpha \vee \tilde{\pi}}\\
	{\sph^4} &
	}
$$
In terms of the diagram (6) it is not hard to verfiy that
\begin{equation}
	\H(f(\alpha \tilde{\pi})) = n^2\cdot \H(\alpha) +
	\H(f(\tilde{\pi}))
\end{equation}
where $\H(g)$ denotes the Hopf invariant of the map $g$.
\medskip

Let $F$ be the homotopy fiber of $\tilde{\pi}$. We will calculate the
cohomology of $F$ using the Serre spectral sequence for the fibration,
\begin{equation}
	\Omega \sph^4 \rightarrow F \rightarrow X
\end{equation}
Recall that $H^\ast(\Omega \sph^4) = \langle z_{3k}, k=0,1,2,\ldots
\rangle$. Consider the diagram of fibrations:
\begin{equation}
\begin{aligned}
\begin{CD}
	\Omega \sph^4 @= \Omega \sph^4 \\
	@VVV @VVV \\
	F @>>> \ast \\
	@VVV @VVV \\
	X @>{\tilde{\pi}}>> \sph^4
\end{CD}
\end{aligned}
\end{equation}
The next proposition is an easy consequence of the naturality of the
Serre spectral sequence with respect to maps of fibrations.

\begin{prop} In the Serre spectral sequence for (7), we have
$$
	d_4 (z_{3k}) = y_4 \otimes z_{3k-3}
$$
where $y_4 \in H^4(X)$ is a generator.\end{prop}

It follows that the classes, $nz_{3k} \in E^{0,3k}_4$, survive to the
next stage. It remains to calculate $d_7(nz_{3k})$. Let $G$ be the
homotopy fiber of $f(\tilde{\pi})$. Using (6) we get a diagram of
fibrations:
\begin{equation}
\begin{aligned}
\begin{CD}
	\Omega \sph^4 @>{\Omega [n]}>> \Omega \sph^4 \\
	@VVV @VVV \\
	F @>>> G \\
	@VVV @VVV \\
	X @>>> \sph^7
\end{CD}
\end{aligned}
\end{equation}
Notice that in the Serre spectral sequence for $\Omega \sph^4
\rightarrow G \rightarrow \sph^7$, we have the identity $d_7(z_6) =
\H(f(\tilde{\pi}))\cdot y_7$ where $y_7 \in H^7(\sph^7)$ is a
generator. Moreover since $(\Omega [n])^\ast (z_6) = n^2 z_6$, using
(10) we have,

\begin{prop} The Hopf invariant, $\H(f(\tilde{\pi}))$, is a multiple of
$n$ i.e.\ $\H(f(\tilde{\pi})) = \lambda n$, and in the Serre spectral
sequence for (8), we have
$$
	d_7(nz_{3k}) = \lambda y_7 \otimes z_{3k-6}
$$
\end{prop}

From Propositions 2.1 and 2.2, we deduce that
\begin{align*}
	H^i(F) &= \Z\cdot y_3 \quad i=3, \quad y_3 =nz_3 \\
		&= \Z_\lambda \cdot y_{7+3k} \quad i=7+3k, \quad
		k=0,1,2,\ldots
\end{align*}
Using the universal coefficients theorem, the homology is:
\begin{align*}
	H_i(F) &= \Z\cdot x_3 \quad i=3,\\
		&= \Z_\lambda \cdot x_{6+3k} \quad i=6+3k,\quad
		k=0,1,2,\ldots
\end{align*}
where $x_3$ and $y_3$ are dual to each other. Now considering the Serre
spectral sequence for the fibration, $F \rightarrow X \rightarrow
\sph^4 $ converging to $H_\ast (X)$ we have:
\begin{align*}
	E^2_{p,q} = E^4_{p,q} &= H_p(\sph^4)\otimes H_q(F) \\
			d_4(x_4) &= n x_3
\end{align*}
where $x_4 \in H_4(\sph^4)$ is a suitably chosen generator. Since
$H_6(X)=0$, the map, $d_4: E^4_{4,3} \rightarrow E^4_{0,6} =
\Z_\lambda$ must be an epimorphism. Hence the class $\lambda
x_4\otimes x_3 \in E^\infty_{4,3}$ represents an orientation $[X] \in H_7(X)$.
\medskip

In the dual picture for the cohomology Serre spectral sequence
converging to $H^\ast(X)$, we have:
\begin{equation}
\begin{aligned}
	E_2^{p,q} = E_4^{p,q} &= H^p(\sph^4)\otimes H^q(F) \\
			d_4(y_3) &= n y_4
\end{aligned}
\end{equation}
where $y_4\in H^4(\sph^4)$ is the class dual to $x_4$. The classes
$y_3$ and $y_4$ are permanent cycles in the spectral sequence
converging to $H^\ast(X,\Z_n)$.
\medskip

From the definition of $[X]$ we have,
\begin{equation}
\langle [y_3],y_4\cap [X] \rangle = \langle [y_3] \cup [y_4],[X]
\rangle \equiv \lambda \!\mod(n)
\end{equation}
As in the previous section, (11) and (12) imply
$$
	\langle \beta^{-1}(y_4),y_4\cap [X] \rangle \equiv \lambda
	\mod(n)
$$
Since the linking form is assumed to be equivalent to a standard form,
it follows that $\lambda \equiv \pm \tau^2 \mod(n)$ for some $\tau \in
(\Z_n)^\ast$. Let $m$ be an integer so that $m\equiv \tau^{-1}
\mod(n)$. Define $\pi: X\rightarrow \sph^4$ as the composite, $\pi =
[m] \comp \tilde{\pi}$. We now have a commutative diagram:
$$
\begin{CD}
	X @>>> \sph^7 \\
	@V{\pi}VV @VV{f(\pi) = [m]\comp f(\tilde{\pi})}V \\
	\sph^4 @>{[n]}>> \sph^4
\end{CD}
$$
Note that $\H(f(\pi))\equiv \pm n\mod(n^2)$. Using (7) we may further assume
$\H(f(\pi)) = \pm n$. It follows now from Proposition 2.2 that
the homotopy fiber of $\pi$ is equivalent to $\sph^3$ since it is a
simply connected homology 3-sphere. We have therefore succeeded in
writing $X$ as the total space of an $\sph^3$ fibration over $\sph^4$
as required.

\section{From fibrations to bundles}

In this section we show that any $\sph^3$ fibration over $\sph^4$ is
equivalent to an $\sph^3$-bundle over $\sph^4$. The argument is fairly
standard, but we outline it for completeness.
\medskip

Let $\xi$ be an $\sph^3$ fibration over $\sph^4$. Restricting $\xi$ to
the hemispheres, $D_+$ and $D_-$, we get trivial fibrations, $\xi_+$
and $\xi_-$, respectively. The fibration $\xi$ may then be obtained by
gluing $\xi_+$ and $\xi_-$ along their common boundary, $\sph^3$, by a
map lying in a homotopy class,
$$
	\sph^3 \rightarrow SG(3)
$$
where $SG(3)$ denotes the monoid of orientation preserving self maps
of $\sph^3$.
\medskip

Recall that (linear) $\sph^3$ bundles are classified by homotopy
classes of maps, $\sph^3 \rightarrow \SO(4)$. The classical $J$
homomorphism identifies $\SO(4)$ with the submonoid of $SG(3)$ of
linear actions. Therefore to show that any $\sph^3$ fibration over
$\sph^4$ is equivalent to an $\sph^3$-bundle over $\sph^4$, it
suffices to show that the map,
$$
	J_\ast : \pi_3 (\SO(4)) \rightarrow \pi_3(SG(3))
$$
is an epimorphism.
\medskip

For a fixed basepoint of $\sph^3$, one has an evaluation map that
evaluates the effect of a self map of $\sph^3$ on the basepoint. It is
easy to see that this map has a section and hence it follows that we
have a map of short exact sequences:
$$
\begin{CD}
0 @>>> \pi_3(\SO(3)) @>>> \pi_3(\SO(4)) @>{\text{ev}_\ast}>>
\pi_3(\sph^3) @>>> 0 \\
@. @VV{J_\ast}V @VV{J_\ast}V @| \\
0 @>>> \pi_3(SG_\ast (3)) @>>> \pi_3(SG(3)) @>{\text{ev}_\ast}>>
\pi_3(\sph^3) @>>> 0
\end{CD}
$$
where $SG_\ast (3)$ are basepoint preserving elements of $SG(3)$.

It suffices to show that $J_\ast : \pi_3(\SO(3)) \rightarrow
\pi_3(SG_\ast (3))$ is an epimorphism. Then the result will follow by
the 5-lemma. Consider the stabilization of $J$:
$$
\begin{CD}
	\pi_3(\SO(3)) @>{J_\ast}>> \pi_3(SG_\ast (3)) \\
	@V{\times 2}VV @VV{\times 2}V \\
	\pi_3 (\SO) @>{J^s_\ast}>> \pi^s_3(S^0)
\end{CD}
$$
It is well known that $\pi_3(SO(3)) = \Z = \pi_3(SO)$. Furthermore, it
is also known that $\pi_3(SG_\ast (3)) = \Z_{12}$ and $\pi^s_3(S^0) =
\Z_{24}$.  By \cite{adams}, $J^s_\ast$ is an epimorphism and hence
$J_\ast$ is an epimorphism as well. This completes the proof of
Theorem 1.

\section{The Berger space}

We briefly recall the construction of the Berger space. Consider
$\R^5$ represented as the space of $3\times 3$ traceless, symmetric
matrices. Then the conjugation action of $\SO(3)$ on this space
affords a (maximal) representation into $\SO(5)$. The quotient space,
$M^7 = \SO(5)/\SO(3)$ may also be written as $\syp(2)/\syp(1)$ for a
maximal embedding of $\syp(1)$ into $\syp(2)$. Berger showed in
\cite{berger} that this space admits a normal homogeneous metric of
positive sectional curvature. In \cite{grovziller} it was shown
that this space cannot be a principal $\sph^3$-bundle over
$\sph^4$. We shall address the question of whether it is equivalent to
an $\sph^3$ fiber bundle over $\sph^4$.
\medskip

The cohomology of this space is well known. However, we outline the
calculation as we will need the setup to compute the linking form. In
terms of the standard maximal tori we
have a commutative diagram:
\begin{equation}
\begin{aligned}
\begin{CD}
	\syp(1) @>{\psi}>> \syp(2) \\
	@AAA @AAA \\
	S^1 @>{\psi_{1,3}}>> S^1 \times S^1
\end{CD}
\end{aligned}
\end{equation}
where $\psi_{1,m}(z) = (z,z^m)$.
\medskip

Let $B\psi : B_{\syp(1)} \rightarrow B_{\syp(2)}$ be the map on the
level of classifying spaces. It follows from (13) that in cohomology
we have:
\begin{equation}
\begin{aligned}
	B\psi^\ast (p_1) &= 10 p_1 \\
	B\psi^\ast (p_2) &= 9 p_1^2
\end{aligned}
\end{equation}
where $H^\ast(B_{\syp(2)}) = \Z[p_1,p_2]$ and $H^\ast(B_{\syp(1)}) =
\Z[p_1]$.
\medskip

The homogeneous space, $M=\syp(2)/\psi(\syp(1))$ is the concrete
description of the Berger space. To calculate its cohomology, consider
the fibration
\begin{equation}
	\syp(2) \rightarrow M \rightarrow B_{\syp(1)}
\end{equation}
We have a pullback diagram:
$$
\begin{CD}
	\syp(2) @= \syp(2) \\
	@VVV @VVV \\
	M @>>> E_{\syp(2)} \\
	@VVV @VVV \\
	B_{\syp(1)} @>{B\psi}>> B_{\syp(2)}
\end{CD}
$$
Recall that $H^\ast (\syp(2)) = \text{E}(y_3,y_7)$, where $y_3$ and
$y_7$ transgress to $p_1$ and $p_2$ respectively in
the Serre spectral sequence for the universal fibration. Using (14)
and the pullback diagram above, we have:

\begin{prop} In the Serre spectral sequence for (15) converging to
$H^\ast(M)$, we have,
\begin{align*}
	d_4(y_3) &= 10 p_1 \\
	d_8(y_7) &= 9 p_1^2
\end{align*}
\end{prop}
It follows immediately from Proposition 4.1 that
\begin{equation}
\begin{aligned}
	H^i(M) &= \Z \quad i=0,7 \\
		&= \Z_{10} \quad i=4
\end{aligned}
\end{equation}
We would like to know whether $M$ is homotopy equivalent to an
$\sph^3$-bundle over $\sph^4$. By Corollary 2, it suffices to compute
the linking form for $M$.
\medskip

Let $S_{1,3} \subset \syp(2)$ denote the $\psi$-image of the standard
maximal torus in $\syp(1)$. We then have a fibration:
\begin{equation}
	\sph^2 \rightarrow \syp(2)/S_{1,3} \rightarrow M
\end{equation}
An easy spectral sequence argument then yields

\begin{prop} $H^\ast(M)$ maps isomorphically onto
$H^\ast(\syp(2)/S_{1,3})$ in degrees 0,4 and 7. The corresponding maps
in homology are isomorphisms as well in degrees 0,3 and 7.\end{prop}

Fix an orientation $[M]\in H_7(M)$. We identify $[M]$ with a class
$[M]$ in $H_7(\syp(2)/S_{1,3})$ using Proposition 4.2. It is clear that
the linking form on $M$ is equivalent to the form,
\begin{align*}
	\alpha: &H^4(\syp(2)/S_{1,3}) \otimes H^4(\syp(2)/S_{1,3})
	\rightarrow \Z_{10} \\
	&x\otimes y \longmapsto \langle \beta^{-1}(x),y\cap [M]\rangle
\end{align*}

It will be easier to calculate $\alpha$ on $\syp(2)/S_{1,3}$ than the
linking form on $M$.
\medskip

Let $\Delta = \syp(1)\times \syp(1) \subset \syp(2)$, be the standard
diagonal embedding of $\syp(1)\times \syp(1)$. Consider the fibration:
\begin{equation}
	\Delta/S_{1,3} \rightarrow \syp(2)/S_{1,3} \rightarrow \syp(2)/\Delta
\end{equation}

The homogeneous spaces $\Delta/S_{1,3}$ and $\syp(2)/\Delta$ can be
identified with the spaces $\sph^2\times \sph^3$ and $\sph^4$
respectively. Hence their cohomologies have the structure of exterior algebras:
$$
H^\ast(\Delta/S_{1,3}) = \text{E}(y_2,y_3), \qquad
H^\ast(\syp(2)/\Delta) = \text{E} (y_4),
$$
where $y_4$ is chosen so that in the Serre spectral sequence for (18),
we have
\begin{equation}
	d_4(y_3) = 10 y_4
\end{equation}

From Proposition 4.2, this is the only non-trivial differential. Since
there are no extension problems, we identify $H^\ast(\syp(2)/S_{1,3})$
with the $E_\infty$ term.  We can do the same for the Serre spectral
sequence in $\Z_{10}$ coefficients, converging to $H^\ast(\syp(2)/S_{1,3},\Z_{10})$. Note that
$y_4 \otimes y_3\in H^7(\syp(2)/ S_{1,3})$ is a generator and hence,
\begin{equation}
\langle [y_3],y_4\cap[M]\rangle = \langle y_4\otimes y_3,[M]\rangle
\equiv \pm 1\!\mod(10)
\end{equation}

As in Section 1, we deduce from (19) and (20) that $\alpha(y_4,y_4) \equiv
\pm 1\mod(10)$. So if $[M]$ is chosen suitably, then the linking form for
the Berger space is equivalent to a standard form.

\appendix
\section{Wilkens' $\beta$ invariant for $\sph^3$-bundles over $\sph^4$}

Equivalence classes of $\sph^3$-bundles over $\sph^4$ with structure
group $\SO(4)$ are in one--one correspondence with $\pi_3(\SO(4))\cong
\Z\oplus \Z$. One can construct generators $\rho,\sigma \in
\pi_3(\SO(4))$ as follows:
$$
	\rho(u)\cdot v = uvu^{-1} , \qquad \sigma(u) = uv
$$
where $u,v$ represent quaternions of norm 1. We now adopt the
convention of \cite{be}. With the above choice of generators, the pair
$(m,n)$ will represent the bundle $\xi_{m,n}$ corresponding to $m\rho
+ n\sigma$. The total space of the associated $\sph^3$-bundle over
$\sph^4$ will be denoted by $M_{m,n}$. The cohomology of $M_{m,n}$ can
be computed without too much difficulty using the Serre spectral
sequence (cf.\ \cite{be});
\begin{align*}
	&H^0(M_{m,n}) \cong H^7(M_{m,n}) = \Z \\
	&H^4(M_{m,n}) = \Z_n \cdot x_4
\end{align*}
where $x_4$ is the pullback of the generator of $H^4(\sph^4)$.
\medskip

In \cite{wilkens}, D. Wilkens studied the class of 2-connected,
7-manifolds. For such manifolds he considered the set of invariants $(H^4(M)^\ast, b, \beta)$. Here $G^\ast$ refers to the torsion part of the group $G$, $b$ is the linking form and $\beta \in H^4(M)$ is the spin characteristic class of the tangent bundle such that $2\beta = p_1$.
\medskip

Since we are interested in $\sph^3$-bundles over $\sph^4$, we restrict 
ourselves to the case when $H^4(M) = \Z_n$, a finite,
cyclic group. Wilkens showed that if $n$ is odd, then the data
$(H^4(M),b,\beta)$ uniquely classifies the manifold $M$ up to oriented
$PL$-homeomorphism type. When $n$ is even, there are at most two
inequivalent manifolds with the same data. For the manifolds $M_{m,n}$
we know that $H^4(M_{m,n}) = \Z_n \cdot x_4$ and the linking form is
standard. Hence the Wilkens invariants for $M_{m,n}$ are
$(H^4(M_{m,n}),\beta)$. From \cite{tamura} it follows that the
characteristic class $\beta$ for $M_{m,n}$ is $\pm 2m x_4$. This is
because the tangent bundle $TM_{m,n}$ is stably equivalent to the
pullback of the vector bundle corresponding to $\xi_{m,n}$.
\medskip

In \cite{jw1}, \cite{jw2}, James and Whitehead have shown that
$M_{m,n}$ is oriented homotopy equivalent to $M_{m',n}$ if and only if
$m\equiv \pm m' \!\mod(n,12)$ (see also \cite{sasao}). When $n =
2\cdot \text{odd}$, the above condition implies that the manifolds $M_{m,n}$ and
$M_{m+\frac{n}{2},n}$ represent distinct oriented homotopy
types. Since they have the same Wilkens data, they realize the two
distinct (oriented) $PL$-homeomorphism types suggested by
\cite[Theorem 2]{wilkens}. Hence in the case $\frac{n}{2}$ is an odd integer, the two distinct possibilities for the data $(H^4,b,\beta)$ are realized by
$\sph^3$-bundles over $\sph^4$.

\begin{proof}[Proof of Theorem 3]: Let $M$ be a simply connected
7-manifold with integral cohomology as in (1) and linking form
equivalent to a standard form. Pick an orientation on $M$ and an
isomorphism, $\psi: \Z_n \rightarrow H^4$ so that $b(\psi(a),\psi(b))
= ab$. Let $x=\psi(1)$; then $b(rx,sx) = rs$.
\medskip

For any manifold with the above cohomology, $\beta \!\mod(2) = w_4$,
the fourth Stiefel--Whitney class of the tangent bundle. An easy
calculation using the Wu classes shows that for such manifolds, the
total Stiefel--Whitney class is trivial and in particular $w_4
=0$. Hence, $\beta = 2m x$, an even class. So the Wilkens data for $M$
is equivalent to the Wilkens data for the manifold $M_{m,n}$. From the
previous discussion, it follows that when $n$ is odd or when $n$ is
twice an odd number, $M$ is $PL$-homeomorphic to $M_{m,n}$.\end{proof}

\begin{rem} Given a manifold that is $PL$-homeomorphic to an $\sph^3$-bundle over $\sph^4$, the obstruction to diffeomorphism is measured by the $\mu$ invariant of Eells and Kuiper (cf.\ \cite{ek}). Computing this invariant requires exhibiting the manifold in question as the boundary of an eight dimensional spin manifold.\end{rem}

\begin{rem} It seems reasonable to expect that the methods outlined here
can be used to prove the analog of Theorem 1 for $\sph^7$-bundles over
$\sph^8$.\end{rem}

\bibliography{bundle}
\bibliographystyle{abbrv}

\end{document}